\documentclass[11pt,a4paper]{article}

\title{Hodge's harmonic $p$-sets and Pontrjagin classes}

\author{Jonathan Fine\relax
\thanks{203 Coldhams Lane, Cambridge, CB1 3HY, England.
\quad E-mail: \texttt{j.fine@pmms.cam.ac.uk}}
}
\date{17 February 1998}

\newcommand\bfR{{\bf R}}
\newcommand\bibrule{\rule{2pc}{0.4pt}}

\begin{document}

\maketitle
\begin{abstract}\noindent
This paper shows how Hodge's theory of harmonic $p$-sets (a discrete
version of his theory of harmonic forms) allows a new approach to be taken
to the problem of providing a combinatorial definition of the Pontrjagin
classes of a compact manifold.  This approach is then related to the
author's definition of flag vectors for hypergraphs, and other objects
constructed out of vertices and cells.
\end{abstract}

\noindent
It is not widely known that Hodge developed a discrete version of his
theory of harmonic differential forms.  The purpose of this paper is to
show how this theory leads to a new approach to the problem of finding a
purely combinatorial description of the Pontrjagin classes of a compact
manifold.  Gelfand and MacPherson \cite{bib.IMG-RDM.CFPC} have already
solved this problem.  Their solution uses oriented matroids, which are
complicated combinatorial objects that can encode `differential
structure'.  The Hodge approach, if successful, is likely to be simpler
and more explicit.  It may also lead to fresh insight into characteristic
classes and related invariants.

In this paragraph, $M$ will be a compact differential manifold of
dimension~$n$, equipped with a  Riemannian metric.  Now suppose that
$\eta$ is any closed $p$-form on $M$, perhaps representing a Pontrjagin
class of $M$.  Hodge's theory of harmonic forms can now be applied.  It
yields a unique $p$-form $\eta_h$, that is equivalent to $\eta$, and
harmonic for the given metric.  Now let $A$ be any $p$-cell of $M$.  The
number
\[
    \eta_A = \int_A\eta_h
\]
is now well-defined, although it does implicitly depend on the choice of a
metric.  Conversely, the system of numbers $\{\eta_A\}$ determines
$\int_B\eta_h$ for any cycle $B$ in $H_p(M,\bfR)$, and so by de~Rahm's
theorem the system $\{\eta_A\}$ also determines the (de~Rham) cohomology
class represented by $\eta$.

Hodge realised that a triangulation is a discrete form of a metric, and
that given a triangulation a discrete version of the harmonic theory can
be developed.  Indeed, he attempted to use the discrete result as the
basis for the proof of the existence of harmonic forms, but he could not
overcome the difficulties inherent in the details of such an approach (see
\cite[p116]{bib.WVDH.TAHI}).

Hodge's theory of harmonic $p$-sets proceeds as follows.  First suppose
that $M$ has been triangulated.  A \emph{$p$-set} is any rule
\[
    A \mapsto \eta_A
\]
that assigns a real number $\eta_A$ to each $p$-cell of the triangulation.
(Previously, we knew $\eta_A$ for all $p$-cells.  Now, we have $\eta_A$
only for those that are in the triangulation.) If $B$ is a $(p+1)$-cell in
the triangulation then the boundary $dB$ of $B$ is a formal sum of
$p$-cells. Say that $\eta$ is \emph{closed} if the sum of $\eta$ evaluated
on $dB$ is zero, for every $(p+1)$-cell $B$.  Each closed $p$-set
determines an element of $H^p(M,\bfR)$, or equivalently a linear function
on $H_p(M,\bfR)$.

Now consider the decomposition of $M$ into cells, that is dual to the
given triangulation.  Each $p$-cell $A$ of the given triangulation will
determine a $(n-p)$-cell $A^*$ of the dual cell decomposition.  If $\eta$
is a $p$-set for the triangulation then the rule
\[
    A^* \mapsto \eta_A
\]
defines a $(n-p)$-set $\eta^*$ for the dual cell decomposition.  Now say
that $\eta$ is \emph{harmonic} (for the given triangulation) if both
$\eta$ and $\eta^*$ are closed.  This is a linear condition on the $p$-set
$\eta$, that depends only on the combinatorial structure of the
triangulation.

Hodge's basic result is that, given a triangulation, each class in
$H^p(M,\bfR)$ has a \emph{unique} representation as a harmonic $p$-set. 
The proof of this result, and the associated definitions, appear on pages
88--92 (\S23.2) and 113--117 (\S28.1) of Hodge's 
book~\cite{bib.WVDH.TAHI}.

Now suppose that $M$ is a compact triangulated manifold, and that $\eta$
is one of its Pontrjagin classes.  Think of $\eta$ as a cohomology class,
and now apply Hodge's theory of harmonic $p$-sets.  Now let $\eta$ be the
harmonic $p$-set that represents the Pontrjagin class.  This construction,
which is part topological and part combinatorial, produces a rule
\[
    A \mapsto \eta_A
\]
that attaches a number $\eta_A$ to each $p$-cell $A$ of the triangulated
manifold.  Although this rule is completely determined by the
combinatorics of $M$, at present its definition requires the intervention
of topology (to define the Pontrjagin class).  To define this rule without
recourse to topology is to provide a purely combinatorial definition of
this Pontrjagin class, via Hodge's theory of harmonic $p$-sets.  This is
the new approach mentioned earlier.

Elsewhere the author has made some definitions that may contribute to the
solution of this problem.  Each triangulated compact $n$-manifold $M$ can
be thought of as a special kind of $(n+1)$-graph $G$.  (The cells of $G$
are the maximal cells of the triangulation.)  In \cite{bib.JF.QTHGFV}, and
more concisely in \cite{bib.JF.SFV}, a flag vector $fG$ is defined for
each such graph.  One can think of $fG$ as an indexed system $\{f_IG\}$ of
integers.  This is not quite appropriate for the present situation.  In
\cite[\S3]{bib.JF.QTHGFV} it is explained how to attach a flag vector to
anything that is built out of vertices and cells.  Thus, for each pair
$A\subset M$ consisting of a $p$-cell $A$ on a triangulated $n$-manifold
$M$, there is a flag vector $f(A\subset M)$ for this pair.  Each component
$f_I$ of this flag vector can, via the rule
\[
    A \mapsto f_I (A\subset M) \>,
\]
be thought of as a $p$-set $\eta_I$, in the sense of Hodge.

Provided the flag vector is in some sense `large enough', each
Hodge-Pontrjagin $p$-set $\eta$ will be a linear combination
\begin{equation}
\label{eqn.eta-sum-aI-etaI}
    \eta = \sum a_I \eta_I
\end{equation}
of the basic flag-vector $p$-sets $\eta_I$.  The main problem can now be
approached as follows.  If $\eta$ is Hodge-Pontrjagin, then $\eta$ is
harmonic ($\eta$ and $\eta^*$ are both closed).  In addition, the class it
determines in $H^p(M,\bfR)$ remains unchanged under subdivision of the
triangulation.  Now determine all $\eta$ of the form
(\ref{eqn.eta-sum-aI-etaI}) that have these properties.  Provided this
combinatorial problem can be solved, all being well one will expect to
find in its solution set the Hodge-Pontrjagin $p$-sets.  One will also
expect to find all products of such.  If this does not exhaust the
solution set, then anything that remains is a candidate to be another
topological (or maybe differential) invariant, of `characteristic class'
type.

The author is most grateful to Laura Anderson and Peter Mani, for helpful
conversations on this problem.  These discussions took place at the recent
Oberwolfach meeting on \emph{Combinatorial Convexity and Algebaic
Geometry}.  The author also thanks the organisers of the meeting for their
invitation, and of course the Oberwolfach Institute itself.

\end{document}